\documentclass[12pt,fullpage]{article}

\usepackage{amsfonts,graphicx,amsmath,xr,epsfig,verbatim,colortbl}

\def\N{\mathbb{N}}
\def\R{\mathbb{R}}

\def\bbP{\mathbb{P}}
\def\bbE{\mathbb{E}}
\def\epsilon{\varepsilon}
\def\tilde{\widetilde}

\def\eps{\epsilon}

\newcommand{\be}{\begin{equation}}
\newcommand{\ee}{\end{equation}}
\newcommand{\baa}{\begin{array}}
\newcommand{\eaa}{\end{array}}
\newcommand{\ba}{\begin{eqnarray}}
\newcommand{\ea}{\end{eqnarray}}

\newtheorem{theo}{\bf Theorem}[section]
\newtheorem{lem}[theo]{\bf Lemma}

\begin{document}
\date{}
\title{\bf{The Harnack inequality for a class of degenerate elliptic operators}}
\author{Fran\c cois Hamel$^{\hbox{\small{ a}}}$$\ $ and Andrej Zlato{\v{s}}$^{\hbox{\small{ b }}}$
\\
\footnotesize{$^{\hbox{a }}$Aix-Marseille Universit\'e \& Institut Universitaire de France}\\
\footnotesize{LATP, Facult\'e des Sciences et Techniques, F-13397 Marseille Cedex 20, France}\\
\footnotesize{$^{\hbox{b }}$Department of Mathematics, University of Wisconsin, Madison, WI 53706, USA}
}
\maketitle

\begin{abstract} 
We prove a Harnack inequality for distributional solutions to a type of degenerate elliptic PDEs in $N$ dimensions.  The differential operators in question are related to the Kolmogorov operator, made up of the Laplacian in the last $N-1$ variables, a first-order term corresponding to a shear flow in the direction of the first variable, and a bounded measurable potential term.  The first-order coefficient is a smooth function of the last $N-1$ variables and its derivatives up to certain order do not vanish simultaneously at any point, making the operators in question hypoelliptic.
\end{abstract}


\section{Introduction}\label{intro}

We prove a Harnack inequality for distributional solutions to the degenerate elliptic PDE
\be\label{5.1}
\Delta_y u+\beta(y)u_x+\gamma(x,y)u=0
\ee
in cylindrical domains in $\R^N$ with axes in the direction of the first variable $x$.  Here $\gamma$ is bounded measurable and $\beta$ is a smooth function   such that the operator
\be\label{5.1a}
L =\sum_{n=1}^{N-1} X_n^2 + X_0 :=  \sum_{n=1}^{N-1} (\partial_{y_n})^2 + \beta(y)\partial_x = \Delta_y + \beta(y)\partial_x 
\ee
satisfies H\" ormander's hypoellipticity condition. That is, vector fields $\{X_n\}_{n=0}^{N-1}$ and their commutators up to certain order span the whole tangent space $\R^N$ at each $(x,y)$.  Moreover, $\beta$ changes sign so that $L$ is not  parabolic, since then the ``elliptic'' Harnack inequality \eqref{5.3} below would not hold in general.  These conditions on $\beta$ are equivalent to hypothesis \eqref{5.2} below  and our result is then as follows:

\begin{theo}\label{T.5.1}
Let $D\subseteq \R^{N-1}$ be open connected and $u:(a,b)\times D\to[0,\infty)$  a bounded distributional solution of \eqref{5.1}
with $\gamma$ bounded measurable and $\beta$ satisfying for some $r\in\mathbb{N}$,
\be\label{5.2}
\beta\in C^{\infty}(D), \ \ {\inf_{D}\,\beta<0<\sup_{D}\,\beta}, \ \ \text{and} \ \
\sum_{0\le |\zeta|\le r} |D^\zeta \beta(y)|>0 \text{ for all $y\in D$}. 
\footnote{For $\zeta=(\zeta_1,\ldots,\zeta_{N-1})\in\N^{N-1}$, we let $|\zeta|=\zeta_1+\cdots+\zeta_{N-1}$ and $D^{\zeta}\beta(y)=\frac{\partial^{|\zeta|}\beta}{\partial y_1^{\zeta_1}\cdots\partial y_{N-1}^{\zeta_{N-1 }}}(y)$.}
\ee
Then for each $[a',b']\subseteq (a,b)$ and bounded  open  $D'$ with $\overline{D'}\subseteq D$,  there is $C>0$,  depending only on $D$, $D'$, $\beta$ and an upper bound on $(a'-a)^{-1},(b-b')^{-1},b'-a'$, and $\|\gamma\|_\infty$,  such that
\be\label{5.3}
\sup_{(a',b')\times D'} u \le C \inf_{(a',b')\times D'} u.
\ee
\end{theo}

{\it Remark.}  We note that $\Delta_y$ could be replaced by any $x$-independent, uniformly elliptic in $y$ operator on $D$, but for the sake of simplicity we  state the theorem with $\Delta_y$ instead.
\smallskip

This result is motivated by its application in our work \cite{hz} on large amplitude $A\to\infty$ asymptotics of traveling fronts in the $x$-direction, and their speeds, for the reaction-advection-diffusion equation
\be\label{eqv}
v_t+A\,\alpha(y)\,v_x=\Delta_{x,y} v+f(v)
\ee
on $\R^{N+1}$, with the first order term representing a shear flow in the $x$-direction and $f$ a non-negative reaction function vanishing at 0 and 1.  The front speeds in question are proved to satisfy $\lim_{A\to\infty}c^*(A\alpha,f)/A= \kappa(\alpha,f)$ for some constant $\kappa(\alpha,f)\ge 0$,  so after substituting the front ansatz $v(t,x,y)=u(x-c^*(A\alpha,f)t,y)$ into \eqref{eqv} and scaling by $A$ in the $x$ variable, one formally recovers \eqref{5.1} in the limit $A\to\infty$, with $\beta(y):=\kappa(\alpha,f)-\alpha(y)$ and $\gamma(x,y):=-f(u(x,y))/u(x,y)$.

The study of hypoelliptic operators of the form 
$$L=\sum_{n=1}^{M} X_n^2 + X_0$$
(where $X_n$ are first order differential operators with smooth coefficients), possibly with an additional potential term,
 has been systematically pursued since H\" ormander's fundamental paper \cite{Hor}.  Although various regularity  and maximum principle results have been obtained soon thereafter (see, e.g., \cite{Bony,FP,Fol,NSW,rs,sv}),  Harnack inequalities and related heat kernel estimates for such operators have initially been proved only in the case when the tangent space at each point is spanned by the fields $\{X_n\}_{n=1}^{M}$ and their commutators, sometimes with $X_0$ either zero or a linear combination of $\{X_n\}_{n=1}^{M}$ \cite{Bony,Jer,JS,KS,KS2}.  

More recently, Harnack inequalities have been obtained without this assumption for certain special classes of operators, not including \eqref{5.1} with general $\beta,\gamma$.  Specifically, some operators with constant and linear coefficients, such as the Kolmogorov operator $L=\partial^2_{yy} + y\partial_x -\partial_t$, were considered in \cite{GL,PP}, and cases of more general coefficients satisfying somewhat rigid structural assumptions (see hypothesis [H.1] in \cite{PP2}) were studied in \cite{KL,PP2} and with a potential term in \cite{PR}.   The domains involved in the obtained inequalities have to depend on the metrics associated to the operators rather than the Euclidian metric, as shows a counter-example to a Harnack inequality in \cite{GL}.  This is related to the need for the sign-changing assumption on $\beta$ here.  We also note that the operators considered in these papers involve the term $\partial_t$ and appropriate ``parabolic-type'' Harnack inequalities are obtained, but corresponding ``elliptic'' inequalities follow from these. 

It was a mild surprise to us that we were not able to find in the literature a sufficiently general result which would include our case \eqref{5.1}.  It appears that Harnack inequalities and heat kernel estimates become much more involved when the field $X_0$ is required for H\" ormander's condition to be satisfied.  One hint in this direction is the fact that the sign-changing hypothesis on $\beta$ is necessary for \eqref{5.3} to hold, so hypoellipticity of $L$ is in itself not a sufficient condition.  

We therefore believe that our method of proof of Theorem \ref{T.5.1} in the next section is itself also a valuable contribution to the problem of quantitative estimates for hypoelliptic operators.  The proof is based on the Feynman-Kac formula for the stochastic process associated with the operator $L$, and  uses the independence of $L$, and thus also of the stochastic process, on $x$.  It is not immediately obvious whether this requirement can be lifted and replaced, for instance, by some assumption on the relation of the stochastic processes associated to $L$ and starting from two different points which can be connected by a path with tangent vector $X_0$ at each point.  We leave this as an open problem.

We thank Luis Caffarelli, Nicola Garofalo, Nicolai Nadirashvili, Brian Street, and Daniel Stroock for useful discussions and pointers to references.  FH is indebted to the Alexander von~Humboldt Foundation for its support.  His work was also supported by the French {\it Agence Nationale de la Recherche} through the project PREFERED.  AZ was supported in part by NSF grants DMS-1113017 and DMS-1056327, and by an Alfred P. Sloan Research Fellowship.  Part of this work was carried out during visits by FH to the Departments of Mathematics of the Universities of Chicago and Wisconsin and by AZ to the Facult\'e des Sciences et Techniques, Aix-Marseille Universit\'e, the hospitality of which is gratefuly acknowledged.

\section{Proof of Theorem \ref{T.5.1}}\label{proof}

Without loss we can assume ${\inf_{D'}\,\beta<0<\sup_{D'}\,\beta}$ and $D'$ connected, after possibly enlarging $D'$.
We will also assume $a=-5$, $a'=0$, $b'=1$, $ b=6$, $D=B_3(0)$, $D'=B_1(0)$, and $\|\gamma\|_\infty\le 1$, with $C$ then only depending on $\beta$, because the general case is analogous.  We also note that \cite[Theorem 18(c)]{rs} and boundedness of $u$ show that $u$ is actually continuous.

We first claim that for each $d>0$ there is  $C_{d,\beta}\ge 1$ such that 
\be\label{5.4}
\sup_{{[0,1]}\times A_d} u \le C_{d,\beta} \inf_{{[0,1]}\times B_1(0)} u,
\ee
with $A_d:=A_d^+\cup A_d^-$ and $A_d^\pm := \{ y\in B_1(0) \,\big|\, \pm\beta(y)> d \}$.  Clearly it suffices to show this for all small enough $d$ such that $A_d^\pm\neq \emptyset$, which we shall assume.

To this end,  note that parabolic regularity theory with $x$ as the time variable, applied on $[-1,5]\times  \{ y\in B_2(0) \,\big|\, -\beta(y)> d/2 \}$, yields
\be\label{5.5}
\sup_{[0,1]\times A_d^-} u \le C_{d,\beta}' \inf_{[2,5]\times A_d^-} u,
\ee
where $C'_{d,\beta}>0$ depends only on $d$ and $\beta$.
Similarly, we obtain
\be\label{5.6}
\sup_{[3,4]\times A_d^+} u \le C_{d,\beta}' \inf_{[-1,2]\times A_d^+} u,
\ee

Next, consider the stochastic process $(X_t^{x,y},Y_t^{x,y})$ starting at $(x,y)\in\R\times B_2(0)$ and satisfying the stochastic differential equation
\[
(dX_t^{x,y},dY_t^{x,y})= (\beta(Y^{x,y}_t)dt,  \sqrt{2}\,dB_t) , \qquad (X_0^{x,y},Y_0^{x,y})=(x,y).
\]
Here $t$ is a new time variable and $B_t$ is a normalized Brownian motion on $\R^{N-1}$ with $B_0=0$ (defined on a probability space $(\Omega,\mathcal{B},\bbP)$). We then have
\be\label{5.8b}
(X_t^{x,y},Y_t^{x,y}) = (X_t^{0,y}+x, \sqrt 2 B_t +y).
\ee
for any $(x,y)\in\R\times B_2(0)$,  in particular, $Y_t^{x,y}$ is independent from $x$.  For any $y\in B_2(0)$ we also define the stopping time
\[
\tau=\tau_y:=\inf \big\{t>0\,\big|\, Y_t^{x,y}\notin B_2(0) \big\}.
\]
If $t\wedge \tau:=\min\{t,\tau\}$, then by the Feynman-Kac formula, $\|\gamma\|_\infty\le 1$, and  the parabolic comparison principle, we have for each $t\ge 0$ and $(x,y)\in\R\times B_2(0)$,
\be \label{5.8}
 e^{-t}  \bbE(u(X_{t\wedge \tau}^{x,y},Y_{t\wedge \tau}^{x,y})) \le u(x,y) \le e^{t} \bbE(u(X_{t\wedge \tau}^{x,y},Y_{t\wedge \tau}^{x,y})).
\ee
(The Feynman-Kac formula is usually stated for $C^2$ functions so we provide a proof of \eqref{5.8} in Lemma \ref{L.5.2} below.)
Here
\be \label{5.8a}
\bbE(u(X_{t\wedge \tau}^{x,y},Y_{t\wedge \tau}^{x,y})) = \int_\Omega u(X_{t\wedge \tau}^{x,y}(\omega),Y_{t\wedge \tau}^{x,y}(\omega)) d\bbP(\omega) = \int_{\R\times \overline{B_2(0)}} u(x',y') d\mu_t^{x,y}(x',y'),
\ee
with the probability measure $\mu_t^{x,y}$ on $\R\times \overline{B_2(0)}$ such that $\mu_t^{x,y}(A) = \bbP((X_{t\wedge \tau}^{x,y},Y_{t\wedge \tau}^{x,y})\in A)$ for  Borel sets $A\subseteq\R\times\overline{B_2(0)}$. Notice that  $\mu_t^{x,y}$ is supported on $[x-\|\beta\|_\infty t, x+\|\beta\|_\infty t] \times \overline{B_2(0)}$ and  $\mu_t^{x,y}(\R \times\partial B_2(0)) = \bbP(\tau_y\le t)$.

By \eqref{5.8b}, translation in $x$ equally translates the $\mu_t^{x,y}$, and the ($x$-independent) measure on $\overline{B_2(0)}$ given by $\nu_t^y(A)=\mu_t^{x,y} (\R\times A)$ is just the law of $\sqrt 2 B_{t\wedge\tau_y} +y$, the Brownian motion on $B_2(0)$ starting at $y$, with stopping time $\tau_y$, and with time scaled by a factor of two.  In particular for each $t>0$ there is $h_t>0$ such that for any $y_1,y_2\in B_1(0)$ and any Borel sets $A_1\subseteq B_1(0)$ and $A_2\subseteq \overline{B_2(0)}$,
\be\label{5.8c}
h_t \nu_t^{y_1}(A_2) \le \nu_t^{y_2}(A_2)\le h_t^{-1} \nu_t^{y_1}(A_2) \qquad \text{and} \qquad h_t |A_1|\le \nu_t^{y_1}(A_1)  \le h_t^{-1}|A_1|.
\ee

From this it follows for $t:=\|\beta\|_\infty^{-1}$ that
\be \label{5.9}
\inf_{[0,1]\times B_1(0)} u \ge C_{d,\beta}'' \inf_{[-1,2]\times A_d^+} u
\ee
with  $C_{d,\beta}'':=e^{-t}h_t \min\{|A_d^+|,|A_d^-|\}$ and, similarly, we obtain
\be \label{5.10}
\inf_{[3,4]\times B_1(0)} u \ge C_{d,\beta}'' \inf_{[2,5]\times A_d^-} u.
\ee
Using \eqref{5.5}, \eqref{5.10}, \eqref{5.6}, and \eqref{5.9} (in that order) yields 
\[
\sup_{[0,1]\times A_d^-} u \le  C_{d,\beta} \inf_{[0,1]\times B_1(0)} u,
\]
with $C_{d,\beta}>0$ depending only on $d$ and $\beta$. An analogous argument gives
\[
\sup_{[0,1]\times A_d^+} u \le  C_{d,\beta} \inf_{[0,1]\times B_1(0)} u,
\]
and \eqref{5.4} follows.  

Next we let $v(x,y):=\int_{-z}^z u(x+s,y)ds$ for some $z\in(0,1/3]$ 
\be\label{5.14}
\sup_{[0,1]\times B_1(0)} v \le \tilde C_{z,\beta} \inf_{[0,1]\times B_1(0)} u
\ee
holds  for some $\tilde C_{z,\beta}\ge 1$.  Indeed, it follows from  \eqref{5.8b}, \eqref{5.8}, \eqref{5.8a} that for each $(x,y)\in\R\times B_2(0)$,
\[
 e^{-t}  \int_{\R\times\overline{B_2(0)}} u(x',y') d\mu_t^{x,y;z}(x',y') \le v(x,y) \le  e^t  \int_{\R\times\overline{B_2(0)}} u(x',y') d\mu_t^{x,y;z}(x',y'),     
 \]
where $\mu_t^{x,y;z}(x',y') = \mu_t^{x,y}(x',y') * ( \chi_{[-z,z]}(x')dx'\delta_0(y'))$.  The above claims about $\mu_t^{x,y}$ and the definition of $\nu_t^y$ imply that 
\[
\mu_t^{x,y;z}(x',y') \le \kappa^{x;z}_{t}(x')\times \nu_t^y(y') \le \sum_{m=-M}^M \mu_t^{x+2mz,y;z}(x',y'),
\]
where $\kappa^{x;z}_{t}$ is the measure on $\R$ with $\kappa^{x;z}_{t}(B)=|B\cap [x-z-\|\beta\|_\infty t,x+z+\|\beta\|_\infty t]|$ for any Borel set $B\subseteq\R$, and $M$ is such that $(2M+1) z \ge 2(z+ \|\beta\|_\infty t)$, for instance, $M:=\lceil 1/2+\|\beta\|_\infty t/z \rceil$.  This and the first claim in  \eqref{5.8c} means that
\be\label{5.16}
v(x,y_1) \le e^{2t} h_t^{-2} \sum_{m=-M}^M  v(x+2mz,y_2)
\ee
for any $x\in\R$, $y_1,y_2\in B_1(0)$ and $t>0$.  

Now we take any $x\in[0,1]$, $y_1\in B_1(0)$, and $y_2\in A_d$ for some fixed $d>0$ such that $A_d^\pm\neq \emptyset$.  Pick $t:= (2\|\beta\|_\infty)^{-1}z$ and $M=1$ to obtain using \eqref{5.16},
\[
v(x,y_1) \le e^{2t} h_t^{-2} \int_{-3z}^{3z}  u(x+s,y_2) ds \le e^{2t} h_t^{-2} \int_{-1}^{2}  u(x',y_2) dx'.
\]
Since \eqref{5.4} and its shifts in $x$ give for $c=-1,0,1$,
\[
\sup_{[c-1,c]\times A_d} u \le C_{d,\beta} \inf_{\{c\}\times A_d} u \le C_{d,\beta} \sup_{[c,c+1]\times A_d} u, 
\]
\[
\sup_{[c-1,c]\times A_d} u \le C_{d,\beta} \inf_{\{c-1\}\times A_d} u \le C_{d,\beta} \sup_{[c-2,c-1]\times A_d} u, 
\]
 we obtain  \eqref{5.4} with $[0,1]$ and $C_{d,\beta}$ replaced by $[-1,2]$ and $C_{d,\beta}^3$. This proves \eqref{5.14}. 
Similarly, \eqref{5.14} with $[-1,0]$ and $[1,2]$ in place of $[0,1]$, together with \eqref{5.4}, yield
\[
\sup_{[-1,2]\times B_1(0)} v \le \tilde{C}_{z,\beta}C_{d,\beta} \inf_{[0,1]\times B_1(0)} u.
\]
In a similar way one can also obtain
\be\label{5.17}
\sup_{[-1,2]\times \overline{B_2(0)}} v \le C_{z,d,\beta} \inf_{[0,1]\times \overline{B_2(0)}} u.
\ee
for some  $C_{z,d,\beta}>0$ (recall that $B_2(0)\subset\subset D=B_3(0)$). 

We will now need to use \eqref{5.2} to finish the proof.  This assumption makes the differential operator on the left-hand side of \eqref{5.1} hypoelliptic in the sense of H\" ormander.  It follows that for $t>0$, the measure $\mu_t^{x,y}$ is absolutely continuous when restricted to $\R\times B_2(0)$ and also to $\R\times \partial B_2(0)$ (as an $(N-1)$-dimensional measure in the latter case),  with densities $p_t(x,y,\cdot,\cdot),q_t(x,y,\cdot,\cdot)\ge 0$ such that
$$p_t(x,y,x',y')=p_t(0,y,x'-x,y'),$$
$$q_t(x,y,x',y')=q_t(0,y,x'-x,y'),$$
and  $p_t,q_t$ are bounded functions when restricted to $y\in B_1(0)$ (with $y'\in B_2(0)$ for $p_t$ and $y'\in \partial B_2(0)$ for $q_t$).  For $p_t$ this follows from the same claim for the corresponding measure $\tilde\mu_t^{x,y}$ on $\R^N$ given by \eqref{5.8a} with $t$ in place of $t\wedge\tau$ and $\beta$ smoothly extended to a periodic function on $\R^{N-1}$ (whose density is smooth in all arguments, \cite[Theorem 3]{IK}).  This is because $\tilde\mu_t^{x,y}(A)\ge \mu_t^{x,y}(A)$ for any Borel set $A\subseteq \R\times B_2(0)$.  

For $q_t$ this would follow from the same claim for the corresponding escape measure $\tilde\mu_{\tau}^{x,y}$ on $\R\times \partial B_2(0)$ given by \eqref{5.8a} with $\tau=\tau_y$ in place of $t\wedge\tau$.  We know of such a result for bounded domains only \cite[Corollary 2.11]{BKS} but since $\mu_t^{x,y}$ is supported on a bounded cylinder, it applies in our case as well.  Specifically, take any $a_- < - \|\beta\|_\infty t$ and $a_+ >  \|\beta\|_\infty t$.  There is a  convex open domain  $G$ with a smooth boundary whose intersection with $[a_-, a_+] \times \R^{N-1}$ is $[a_-, a_+] \times B_2(0)$, and the intersection with $[(-\infty, a_-)\cup(a_+,\infty)] \times \R^{N-1}$ are two smooth ``slanted'' conical caps $G_\pm\subseteq \R\times B_2(0)$ over the $(N-1)$-dimensional balls $\{a_\pm\}\times B_2(0)$ with the two (rounded) tips at points with $y'$ coordinates $y'_{\pm}$ such that $\pm \beta(y'_{\pm})>0$ and sufficiently long so that for any $(x',y')\in\partial G_\pm\cap\partial G$, the unit outer normal vector $n(x',y')$ to $\partial G_\pm$ at $(x',y')$ satisfies
$$|n(x',y')\cdot (1,0,\cdots,0)| \le \tfrac 12 (\|\beta\|_\infty^{-1}+1)\ \hbox{ whenever  }\pm\beta(y')\le 0.$$
Then $G$ satisfies the hypotheses of \cite[Corollary 2.11]{BKS} (it satisfies the escape condition and all points of $\partial G$ are $\tau'$-regular).  It follows that the escape measure $\tilde\mu_{\tau}^{x,y}$ has a density $\tilde q_\tau(x,y,\cdot,\cdot)$ which is a continuous function of $(x,y,x',y')\in G\times\partial^* G$, where $\partial^* G$ is the set of ``good'' points of $\partial G$, that is,  all $(x',y')\in\partial G$ except of the two cone tips, where $n(x',y')=(\pm 1,0,\cdots,0)$.  Thus $\tilde q_\tau$ is bounded on $S:=\{0\}\times B_1(0)\times(a_-,a_+)\times\partial B_2(0)$.  Since $\{X^{0,y}_s\}_{s\le t\wedge \tau}$ almost surely stays in $(a_-,a_+)$, we obtain $q_t\le \tilde q_\tau$ on $S$ and $q_t=0$ on $[\{0\}\times B_1(0)\times\R \times\partial B_2(0)] \setminus S$.   Finally, $q_t(x,y,x',y')=q_t(0,y,x'-x,y')$ shows that $q_t$ is bounded on $\R\times B_1(0)\times \R \times\partial B_2(0)$.

Let $d>0$ be such that $A_d^\pm\neq \emptyset$, let $z:=1/3$, $t:=\|\beta\|_\infty^{-1}$, and 
\[
C_t:=\max\{ \sup_{\R\times B_1(0)\times \R \times B_2(0)} p_t,  \sup_{\R\times B_1(0)\times \R \times\partial B_2(0)} q_t\} <\infty.
\]
 Then $p_t(x,y,x',y'),q_t(x,y,x',y')\le C_t \chi_{[x-1,x+1]}(x')$ because the measure  $\mu_t^{x,y}$ is supported on $[x-1, x+1] \times \overline{B_2(0)}$, so we obtain from \eqref{5.8} and \eqref{5.8a}
\begin{align*}
\sup_{[0,1]\times B_1(0)} u & \le C_t e^t\int_{[-1,2]\times B_2(0)} u(x',y') dx'dy' + C_t e^t\int_{[-1,2]\times \partial B_2(0)} u(x',y') dx'dy' 
\\ & \le 10C_te^t \sup_{[-1,2]\times \overline{B_2(0)}} v 
\\ & \le 10C_t C_{z,d,\beta} e^t \inf_{[0,1]\times B_1(0)} u
\end{align*}
by using $[-1,2]=[-1,-1/3]\cup[-1/3,1/3]\cup[1/3,1]\cup[1,5/3]\cup[4/3,2]$ and \eqref{5.17}.  This is \eqref{5.3}, so the theorem will be proved once we establish \eqref{5.8}. 

\begin{lem} \label{L.5.2}
If $u,\beta,\gamma,X^{x,y}_t,Y^{x,y}_t,\tau_y$ are as in the proof of Theorem \ref{T.5.1} (in particular, $\|\gamma\|_\infty\le 1$), then \eqref{5.8} holds for $(x,y)\in \R\times B_2(0)$.
\end{lem}

\noindent{\bf{Proof.}}
Let $Z^{x,y}_t=t$ so that $dZ^{x,y}_t=dt$ and $K:=\Delta_y+\beta(y)\partial_x+\partial_z$ is the generator of the process $(X^{x,y}_t,Y^{x,y}_t,Z^{x,y}_t)$.  If we let $v(x,y,z):=e^{z}u(x,y)$, then $Kv\ge 0$ on $\R\times B_3(0)\times\R$ in the sense of distributions, that is, 
\[
\int_{\R\times B_3(0)\times\R} v K^*\phi\, dxdydz \ge 0
\]
for any $\phi\in C^\infty_0(\R\times B_3(0)\times\R)$, with $K^*:= \Delta_y-\beta(y)\partial_x-\partial_z$ the adjoint of $K$.

For any $\eps>0$ let $\delta_\eps\in(0,1/2\sqrt{N-1})$ be such that $|\beta(y)-\beta(y')| \le \eps^2$ whenever $y,y'\in B_{5/2}(0)$ and  $|y-y'|\le \sqrt{N-1}\,\delta_\eps$.   Let $\eta:\R\to [0,1]$ be a smooth non-negative function supported in $[-1,1]$, with $\int_{-1}^1\eta(x')dx'=1$ and  $\|\eta'\|_\infty\le 2$.  For $\eps>0$ define the mollifier
\[
\eta^{\eps}(x,y,z) :=
\eps^{-2}\delta_\eps^{1-N} \eta \left( \displaystyle{\frac x\eps} \right) \eta \left( \displaystyle{\frac z\eps} \right)  \displaystyle{\mathop{\prod}_{n=1}^{N-1}} \eta \left( \displaystyle{\frac {y_n}{\delta_\eps}} \right),
\]
and let $v^{\eps}:= v*\eta^\eps$ and $\phi^{\eps;x,y,z}(x',y',z'):=\eta^\eps(x-x',y-y',z-z')$. 
For $\eps\in(0,1)$ the smooth function $v^\eps$ then satisfies on $\R\times B_2(0)\times\R$
\begin{align*}
(Kv^{\eps}) (x,y,z) = & \int_{\R\times B_3(0)\times\R} vK^*\phi^{\eps;x,y,z} \, dx' dy' dz'  
\\ + & \int_{\R\times B_3(0)\times\R} v(x',y',z')[\beta(y')-\beta(y)] \phi^{\eps;x,y,z}_{x'}(x',y',z') \, dx' dy' dz'.
\end{align*}
The first integral is non-negative.  The integrand in the second vanishes when $|x'-x|>\eps$ or $|y'_n-y_n|> \delta_\eps$ for some $n$ or $|z'-z|>\eps$, and $|\phi^{\eps;x,y,z}_{x'}(x',y',z')|\le 2 \eps^{-3}\delta_\eps^{1-N}$, so we have
\[
(Kv^{\eps}) (x,y,z) \ge - 2^{N+2} \eps e^{z+\eps} \|u\|_\infty.
\]

We next apply Dynkin's formula \cite[Theorem 7.4.1]{Oks} to the smooth function $v^\eps$, the process $(X^{x,y}_t,Y^{x,y}_t,Z^{x,y}_t)$, and stopping time $t\wedge \tau$ (with $\tau=\tau_y$), to obtain
\begin{align*}
\bbE \left[ v^\eps(X^{x,y}_{t\wedge \tau},Y^{x,y}_{t\wedge \tau},Z^{x,y}_{t\wedge \tau}) \right] & = v^\eps(x,y,0)+ \bbE\left[ \int_0^{t\wedge \tau} (Kv^\eps) (X^{x,y}_s,Y^{x,y}_s,Z^{x,y}_s) ds \right] 
\\ & \ge v^\eps(x,y,0) - 2^{N+2}  \eps e^{t+\eps} \|u\|_\infty t.
\end{align*}
Since $v^\eps\to v$ uniformly on $[x-\|\beta\|_\infty t, x+\|\beta\|_\infty t]\times \overline{B_2(0)}\times [0,t]$ as $\eps\to 0$ (by continuity of $v$) and $Z^{x,y}_{t\wedge \tau}\le t$, it follows that
\[
e^t \bbE \left[ u(X^{x,y}_{t\wedge \tau},Y^{x,y}_{t\wedge \tau}) \right] \ge \bbE \left[ v(X^{x,y}_{t\wedge \tau},Y^{x,y}_{t\wedge \tau},Z^{x,y}_{t\wedge \tau}) \right]\ge u(x,y).
\]
This is the second inequality in \eqref{5.8}.  The first inequality is obtained in the same way, this time with $v(x,y,z):=e^{-z}u(x,y)$, so that $Kv\le 0$ on $\R\times B_3(0)\times\R$ and
\[
(Kv^{\eps}) (x,y,z) \le  2^{N+2} \eps e^{-z+\eps} \|u\|_\infty.
\]
\hfill$\Box$


\end{document}